\documentclass[12pt]{article}  
\usepackage{amsthm,amsmath,amsfonts,epsfig,amssymb,latexsym}
\usepackage{mathrsfs}
\usepackage[T1]{fontenc}
\usepackage{graphicx}
\usepackage{enumerate}
\usepackage{xy,xypic}
\usepackage{graphics}
\usepackage{footnote}
\usepackage{scalefnt}
\usepackage{tikz}
\usetikzlibrary{shapes,arrows,fit,calc,positioning}
\usetikzlibrary{matrix}
\xyoption{all}
\title{{\large An Example in Complete Intersections and an Erratum}
} 
\author{   
 Satya Mandal 
 \\ 
{\small University of Kansas, Lawrence, Kansas 66045;}
{\small {\it  mandal@ku.edu 
  }}\\
 }  
\begin{document}
\renewcommand{\baselinestretch}{1.255}
\setlength{\parskip}{1ex plus0.5ex}
\date{31 January 2017}
\newtheorem{theorem}{Theorem}[section]
\newtheorem{proposition}[theorem]{Proposition}
\newtheorem{lemma}[theorem]{Lemma}
\newtheorem{corollary}[theorem]{Corollary}
\newtheorem{construction}[theorem]{Construction}
\newtheorem{notations}[theorem]{Notations}
\newtheorem{question}[theorem]{Question}
\newtheorem{example}[theorem]{Example}
\newtheorem{definition}[theorem]{Definition} 
\newtheorem{clarification}[theorem]{Clarification} 
\newtheorem{Myproof}[theorem]{Proof} 
\newtheorem{remarks}[theorem]{Remarks}
\newtheorem{remark}[theorem]{Remark}

\newtheorem{conjecture}[theorem]{Conjecture}
\newtheorem{openProblem}[theorem]{Open Problem}

\newcommand{\iso}{\stackrel{\sim}{\longrightarrow}}
\newcommand{\sur}{\twoheadrightarrow}
\newcommand{\bD}{\begin{definition}}
\newcommand{\eD}{\end{definition}}
\newcommand{\bP}{\begin{proposition}}
\newcommand{\eP}{\end{proposition}}
\newcommand{\bL}{\begin{lemma}}
\newcommand{\eL}{\end{lemma}}
\newcommand{\bT}{\begin{theorem}}
\newcommand{\eT}{\end{theorem}}
\newcommand{\bC}{\begin{corollary}}
\newcommand{\eC}{\end{corollary}} 
\newcommand{\eop}{\hfill \rule{2mm}{2mm}}
\newcommand{\pf}{\noindent{\bf Proof.~}}
\newcommand{\PD}{\text{proj} \dim}
\newcommand{\lra}{\longrightarrow}
\newcommand{\hra}{\hookrightarrow}
\newcommand{\llra}{\longleftrightarrow}
\newcommand{\Lra}{\Longrightarrow}
\newcommand{\Llra}{\Longleftrightarrow}
\newcommand{\bE}{\begin{enumerate}}
\newcommand{\eE}{\end{enumerate}}
\newcommand{\pic}{The proof is complete.}
\def\spec#1{\mathrm{Spec}\left(#1\right)}
\def\m{\mathfrak {m}}
\def\Sch{\underline{\mathrm{Sch}}}
\def\Sets{\underline{\mathrm{Sets}}}

\def\CA{\mathcal {A}}
\def\CB{\mathcal {B}}
\def\CP{\mathcal {P}}
\def\CC{\mathcal {C}}
\def\CD{\mathcal {D}}
\def\CE{\mathcal {E}}
\def\CF{\mathcal {F}}
\def\CE{\mathcal {E}}
\def\CG{\mathcal {G}}
\def\CH{\mathcal {H}}
\def\CI{\mathcal {I}}
\def\CJ{\mathcal {J}}
\def\CK{\mathcal {K}}
\def\CL{\mathcal {L}}
\def\CM{\mathcal {M}}
\def\CN{\mathcal {N}}
\def\CO{\mathcal {O}}
\def\CP{\mathcal {P}}
\def\CQ{\mathcal {Q}}
\def\CR{\mathcal {R}}
\def\CS{\mathcal {S}}
\def\CT{\mathcal {T}}
\def\CU{\mathcal {U}}
\def\CV{\mathcal {V}}
\def\CW{\mathcal {W}}
\def\CX{\mathcal {X}}
\def\CY{\mathcal {Y}}
\def\CZ{\mathcal {Z}}

\newcommand{\smallcirc}[1]{\scalebox{#1}{$\circ$}}
\def\BA{\mathbb {A}}
\def\BB{\mathbb {B}}
\def\BC{\mathbb {C}}
\def\BD{\mathbb {D}}
\def\BE{\mathbb {E}}
\def\BF{\mathbb {F}}
\def\BG{\mathbb {G}}
\def\BH{\mathbb {H}}
\def\BI{\mathbb {I}}
\def\BJ{\mathbb {J}}
\def\BK{\mathbb {K}}
\def\BL{\mathbb {L}}
\def\BM{\mathbb {M}}
\def\BN{\mathbb {N}}
\def\BO{\mathbb {O}}
\def\BP{\mathbb {P}}
\def\BQ{\mathbb {Q}}
\def\BR{\mathbb {R}}
\def\BS{\mathbb {S}}
\def\BT{\mathbb {T}}
\def\BU{\mathbb {U}}
\def\BV{\mathbb {V}}
\def\BW{\mathbb {W}}
\def\BX{\mathbb {X}}
\def\BY{\mathbb {Y}}
\def\BZ{\mathbb {Z}}

\newcommand{\TCP}{\textcolor{purple}}
\newcommand{\TCM}{\textcolor{magenta}}
\newcommand{\TCR}{\textcolor{red}}
\newcommand{\TCB}{\textcolor{blue}}
\newcommand{\TCG}{\textcolor{green}}

\def\SA{\mathscr {A}}
\def\SB{\mathscr {B}}
\def\SC{\mathscr {C}}
\def\SD{\mathscr {D}}
\def\SE{\mathscr {E}}
\def\SF{\mathscr {F}}
\def\SG{\mathscr {G}}
\def\SH{\mathscr {H}}
\def\SI{\mathscr {I}}
\def\SJ{\mathscr {J}}
\def\SK{\mathscr {K}}
\def\SL{\mathscr {L}}
\def\SN{\mathscr {N}}
\def\SO{\mathscr {O}}
\def\SP{\mathscr {P}}
\def\SQ{\mathscr {Q}}
\def\SR{\mathscr {R}}
\def\SS{\mathscr {S}}
\def\ST{\mathscr {T}}
\def\SU{\mathscr {U}}
\def\SV{\mathscr {V}}
\def\SW{\mathscr {W}}
\def\SX{\mathscr {X}}
\def\SY{\mathscr {Y}}
\def\SZ{\mathscr {Z}}

\def\bfA{{\bf A}}
\def\bfB{{\bf B}} 
\def\bfC{{\bf C}} 
\def\bfD{{\bf D}} 
\def\bfE{{\bf E}} 
\def\bfF{{\bf F}} 
\def\bfG{{\bf G}} 
\def\bfH{{\bf H}} 
\def\bfI{{\bf I}} 
\def\bfJ{{\bf J}} 
\def\bfK{{\bf K}} 
\def\bfL{{\bf L}} 
\def\bfM{{\bf M}} 
\def\bfN{{\bf N}} 
\def\bfO{{\bf O}} 
\def\bfP{{\bf P}} 
\def\bfQ{{\bf Q}} 
\def\bfR{{\bf R}} 
\def\bfS{{\bf S}} 
\def\bfT{{\bf T}} 
\def\bfU{{\bf U}} 
\def\bfV{{\bf V}} 
\def\bfW{{\bf W}} 
\def\bfX{{\bf X}} 
\def\bfY{{\bf Y}} 
\def\bfZ{{\bf Z}} 

\maketitle

\section{Introduction}

The following is a version of the complete intersection conjecture of M. P. Murthy (\cite{M}, \cite[pp 85]{M1}).

\begin{conjecture} \label{MurthyConj} 
Suppose $A=k[X_1, X_2, \ldots, X_n]$ is a polynomial ring over a field $k$.
Then, for any ideal $I$ in $A$, $\mu(I)=\mu(I/I^2)$, {\rm where $\mu$ denotes the minimal number of generators}.
\end{conjecture}


Recently, a solution of this conjecture (\ref{MurthyConj}) was claimed in \cite{F}, when $k$ is an infinite  perfect field and $1/2\in k$.
Subsequently, using Popescu's Desingularization Theorem (\cite{P, Sw1}), the claim was strengthened in \cite{M0}, removing the perfectness condition. 
In deed,  the main claims in \cite{F}, were further strengthened in \cite{M0}, for ideals $I$ in polynomial rings $A[X]$, over
 regular rings $A$ containing an infinite field $k$, with $1/2\in k$, and $I$ containing a monic polynomial. 
 It was latter established (see \cite{MM}) that the 
 methods in \cite{F, M0}, would also work when $k$ is finite, with $1/2\in k$.

The  claimed proof of Conjecture \ref{MurthyConj} in \cite{F} is a consequence of a stronger claim that, for integers $n\geq 2$,
 any set of $n$-generators of $I/I^2$ lifts to a set of $n$-generators of $I$. This question of liftability of generators of $I/I^2$ was considered in
 \cite{MMu} and a counter example was given \cite[Example 2.4]{MMu}, when $n=2$. 
 In next section  \ref{mainExamp}, we develop a   larger class of examples,
 for all integers $n\geq 2$.

Due to the existence of such counter examples, some clarifications are needed regarding the published claims in \cite{F, M0}.
We provide the same in  Section \ref{erratum}.
We underscore that, there is no logical error in the methods in \cite{M0}, barring the use of the claimed results in  \cite{F}. 


\section{The Example}\label{mainExamp}
The following two examples were worked out in collaboration with M. P. Murthy.

\begin{example}\label{exam} 
Let $n\geq 3$ be any integer, and
$
A=k[X_1, \ldots,X_n; Y_1, \ldots, Y_n]
$
 be a polynomial ring over any field $k$. Let 
$f=\sum_{i=1}^nX_iY_i-1\in A$
 and $I=Af$.
Write $\overline{A}=\frac{A}{(f)}$. For elements in $A$ {\rm (respectively, in $I$)}, the
 images in $\overline{A}$  {\rm (respectively, in $\frac{I}{I^2}$)}
  will be denoted by "overline". Then, 
$$
 \overline{X_1f}, \overline{X_2f}, ~~\ldots,~~ \overline{X_nf}
\quad {\rm generates} \quad \frac{I}{I^2}.
$$
This set of generators of $\frac{I}{I^2}$ would not lift to a set of generators of $I$.
\end{example}

\pf As in \cite{MMu}, we have the commutative diagram 
$$
\diagram 
A\ar[r]^f_{\sim}\ar@{->>}[d]& I\ar@{->>}[d]\\
\frac{A}{(f)} \ar[r]_{\overline{f}}^{\sim} &  \frac{I}{I^2}\\
\enddiagram
$$
Suppose $\overline{X_1f}, \overline{X_2f}, ~~\ldots,~~ \overline{X_nf}$ lifts to a set of generators of $I$.
Then, by the diagram above,  
the unimodular row 
$$
\left(\overline{X_1}, \overline{X_2}, ~~\ldots,~~ \overline{X_n}\right)
\qquad \qquad {\rm of} \quad \overline{A}
$$
lifts to a unimodular row 
$$
(F_1, F_2,  \ldots, F_n) \qquad \qquad {\rm of} \quad A.
$$
 Since projective $A$-modules are free, there is a matrix 
$\sigma\in GL_n(A)$, whose first row is $(F_1, \ldots, F_n)$. Therefore, 
$\left(\overline{X_1}, \overline{X_2}, ~~\ldots,~~ \overline{X_n}\right)$ is the first row of the image of $\sigma$ in 
$GL_n\left(\overline{A}\right)$. 
So, the projective $\overline{A}$-module defined by $\left(\overline{X_1}, \overline{X_2}, ~~\ldots,~~ \overline{X_n}\right)$ is free.
This is impossible, 
by the Theorem of N. Mohan Kumar and Madhav V. Nori 
(see \cite[Theorem 17.1]{Sw1}).
\pic $\eop$

\begin{example}\label{SphereE}
Let $A=\BR[X_0, X_1, \ldots, X_n]$ be a polynomial ring over the field of real numbers. Let $f=\sum_{i=0}^n X_i^2-1\in \BR$ and $I=Af$.
Assume, $n\neq 0, 1, 3, 7$.
Then, $X_0f, X_1f, \ldots, X_nf$ induce a set of generators for $I/I^2$, which would not lift to a set of generators of $I$.
\end{example}
\pf Same as the proof of (\ref{exam}), while we use the fact that tangent bundles over real $n$-spheres ($n\neq 0, 1, 3, 7$)
are nontrivial (see \cite[Theorem 2.3]{Sw1}). $\eop$

\begin{remark}{\rm 
Note $I=Af$ in (\ref{exam}) is a  principal ideal. So, Examples \ref{exam}, \ref{SphereE}, do not provide a counter example of the Complete Intersection
 Conjecture \ref{MurthyConj}.
}
\end{remark}


\section{Erratum} \label{erratum}
The following list provides some clarifications regarding the inconsistencies in the literature \cite{F, M0}, at this time.
\bE
\item  It was communicated by Mrinal K. Das that the proof of 
\cite[Lemma 3.2.3]{F} is not convincing. 
This may be the likely cause of all the inconsistencies in  \cite{F, M0}, under discussion.
In deed, this creates an incompleteness in the proof of the key result   \cite[Theorem 3.2.7]{F}.
\item Theorem 3.2.8 in \cite{F} does not have a valid proof in the literature. For a polynomial ring $A=k[X_1, \ldots, X_n]$, 
by \cite[Proposition 4.1]{M0}, $Q_{2n}(A)$ is a singleton.  Therefore,  \ref{exam}, \ref{SphereE} would be a counter example to  \cite[Theorem 3.2.8]{F}.
\item The claimed proof of  \cite[Theorem 3.2.9]{F} is not valid, since it uses \cite[Theorem 3.2.8]{F}. Therefore, the Complete Intersection Conjecture \ref{MurthyConj}
is still open and the best result on this conjecture, at this time, remains those in \cite{Mk} and \cite{M2}.
\item There is no logical error in \cite{M0}. However, since  the main results in \cite{M0} depends on the validity of the same in \cite{F}, they do 
not have any valid proofs, at this time. In particular, 
\bE
\item The main results \cite[Theorems 3.8, 4.2, 4.3]{M0} do not have  valid proofs, at this time.
\item The claimed proofs of Abhyankar's epimorphism conjecture \cite[Thoerem 4.5, 4.6]{M0}, are not valid, since they 
are routine consequence of \cite[Theorems 4.2, 4.3, 4.4]{M0}.
\eE
\item The results in \cite{M0} that are not dependent on results in \cite{F} are valid. In particular,
\bE
\item The \cite[Propositions 4.1]{M0}, on triviality of homotopy obstructions, for ideals containing a monic polynomial, remain  valid. 
\item Results in \cite[Section 5]{M0}, on the alternate description of the Homotopy Obstruction set $Q_{2n}(A)$, remain  valid.
\eE
\eE

\begin{thebibliography}{200}


\bibitem[F]{F}  Fasel, Jean 
Fasel, Jean On the number of generators of ideals in polynomial rings. {\it Ann. of Math. (2)} 184 (2016), no. 1, 315-331.; arXiv:1507.05734 






\bibitem[MM]{MM}
Mandal Satya, Bibekananda Mishra 
The Homotopy Program in Complete Intersections, arXiv:1610.07495

\bibitem[M0]{M0}
 Mandal, Satya On the complete intersection conjecture of Murthy. {\it J. Algebra} 458 (2016), 156?170.
 
\bibitem[M1]{M1}
Mandal, Satya Projective modules and complete intersections. Lecture Notes in Mathematics, 1672. {\it Springer-Verlag, Berlin}, 1997.

\bibitem[M2]{M2} Mandal, Satya On efficient generation of ideals. Invent. Math. 75 (1984), no. 1, 59-67.
 
 \bibitem[MMu]{MMu}
  Mandal, Satya; Pavaman Murthy, M. Ideals as sections of projective modules. {\it J. Ramanujan Math. Soc.} 13 (1998), no. 1, 51?62.
 
 
 
 \bibitem[Mk]{Mk}
 Kumar, N. Mohan On two conjectures about polynomial rings. {\it Invent. Math}. 46 (1978), no. 3, 225-236.
 
\bibitem[M]{M} Murthy, M. Pavaman Complete intersections. {\it  Conference on Commutative Algebra-1975}
  (Queen's Univ., Kingston, Ont., 1975), pp. 196-211. Queen's Papers on Pure and Applied Math., No. 42,  {\it Queen's Univ., Kingston, Ont.}, 1975.
 
 \bibitem[P]{P}
 Popescu, Dorin Letter to the editor: "General N\'{e}ron desingularization and approximation'' {\it Nagoya Math. J}. 118 (1990), 45-53.;
 
 \bibitem[S]{S}
 Stavrova, A. Homotopy invariance of non-stable $K_1$-functors. {\it J. $K$-Theory} 13 (2014), no. 2, 199-248.
 
\bibitem[Sw1]{Sw1}
 Swan, Richard G. N\'{e}ron-Popescu desingularization. {\it Algebra and geometry (Taipei, 1995)},
  135-192, Lect. Algebra Geom., 2, {\it Int. Press, Cambridge, MA}, 1998.
 

  
 

\end{thebibliography}
\end{document}